
%

\input amstex
\loadbold
\documentstyle{amsppt}
\NoBlackBoxes

\pagewidth{32pc}
\pageheight{44pc}
\magnification=\magstep1

\def\QQ{\Bbb{Q}}

\def\PP{\Bbb{P}}
\def\OO{\Cal{O}}

\def\QSec{\operatorname{QSec}}
\def\Rat{\operatorname{Rat}}
\def\Hilb{\operatorname{Hilb}}
\def\Var{\operatorname{Var}}
\def\id{\operatorname{id}}
\def\rank{\operatorname{rank}}
\def\codim{\operatorname{codim}}
\def\ratmap{\dashrightarrow}
\def\rest#1#2{\left.{#1}\right\vert_{{#2}}}
\def\QED{{\unskip\nobreak\hfil\penalty50\quad\null\nobreak\hfil
{$\square$}\parfillskip0pt\finalhyphendemerits0\par\medskip}}

\topmatter
\title
Remarks on S. Lang's conjecture \\
over function fields
\endtitle
\rightheadtext{}
\author Atsushi Moriwaki \endauthor
\leftheadtext{}
\address
Department of Mathematics, Faculty of Science,
Kyoto University, Kyoto, 606-01, Japan
\endaddress
\curraddr
Department of Mathematics, University of California,
Los Angeles, 405 Hilgard Avenue, Los Angeles, California 90024, USA
\endcurraddr
\email moriwaki\@math.ucla.edu \endemail
\date November, 1994, 1st Draft \enddate
\abstract
In this short note, we will show the following weak evidence
of S. Lang conjecture over function fields.
Let $f : X \to Y$ be a projective and surjective morphism of algebraic
varieties
over an algebraically closed field $k$ of characteristic zero,
whose generic fiber is geometrically irreducible and of general type.
If $f$ is not birationally trivial,
then there are countably many proper closed varieties
$\{ Z_i \}$ of $X$ such that
every quasi-section of $f$ is contained in $\bigcup_i Z_i$.
\endabstract
\endtopmatter

\document

\head \S0. Introduction
\endhead

In \cite{La}, S. Lang conjectured that
if $X$ is a projective variety over a number field $K$ and
$X$ is of general type, then there is a proper subscheme $Z$ of $X$
with $X(K) \subset Z(K)$.
We can consider an analogue over function fields.
In this case, we must avoid a birationally trivial family, i.e.
an algebraic family $f : X \to Y$ of algebraic varieties which is
birationally equivalent to a product $W \times Y$ over $Y$.

\definition{Conjecture A}(Analogue of S. Lang's conjecture over function
fields)
\quad
Let $f : X \to Y$ be a projective and surjective morphism of algebraic
varieties
over an algebraically closed field $k$,
whose generic fiber is geometrically irreducible and of general type.
If $f$ is not birationally trivial,
then there are a proper subscheme $Z$ of $X$
such that every quasi-section of $f$ is contained in $Z$.
\enddefinition

When $\dim f = 1$, Conjecture~A is known as Mordell's conjecture over
function fields and was proved by a lot of authors in any characteristic.
However, Conjecture~A does not hold in this naive form
if $\dim f \geq 2$ and the characteristic is positive.
For, roughly speaking, in this case,
there are a family of unirational varieties of general type (For more details,
see Remark~4.2). On the other hand, if the characteristic is zero,
from all we know, it is still an open problem.
We only know it is true if the cotangent bundle of the generic fiber is
ample (cf. \cite{No} and \cite{Mo}).
In this short note, we will prove the following weak evidence of
the above Conjecture A in characteristic zero.

\proclaim{Proposition B}
Let $f : X \to Y$ be a projective and surjective morphism of algebraic
varieties
over an algebraically closed field $k$ of characteristic zero,
whose generic fiber is geometrically irreducible and of general type.
If $f$ is not birationally trivial,
then there are countably many proper closed varieties
$\{ Z_i \}$ of $X$ such that
every quasi-section of $f$ is contained in $\bigcup_i Z_i$.
\endproclaim

Our basic tool in this note is the following criterion of
birational splitting (cf. Corollary~1.2).
\block{\sl
Let $f : X \ratmap Y$ be a dominant rational map,
whose generic fiber is geometrically irreducible and
of general type. Assume that there exists a variety $T$
and a dominant rational map
$\phi : T \times Y \ratmap X$.
Then $X$ is a birationally equivalent to a product $W \times Y$.}
\endblock
This is a problem raised in Historical appendix of \cite{La} and
essentially was solved by K. Maehara \cite{Ma} earlier than
\cite{La}.

\head \S1. Birational splitting
\endhead

In this section, we will consider a criterion for birational splitting.
It was essentially due to K. Maehara \cite{Ma}.
This is a very important tool for diophantine geometry over
function fields. So we will re-prove it.
Let's us begin with the following lemma, which is an easy application
of weak positivity of direct images of $n$-th relative canonical bundles.

\proclaim{Lemma 1.1}
Let $f : X \to Y$ be a surjective morphism
of smooth projective varieties over an algebraically closed field $k$
of characteristic zero.
If there are a projective smooth algebraic variety $T$ over $k$
and a dominant rational
map $\phi : T \times_k Y \ratmap X$ over $Y$, then
the double dual $f_*(\omega_{X/Y}^n)^{\vee\vee}$ of $f_*(\omega_{X/Y}^n)$
is a free $\OO_Y$-sheaf for all $n \geq 0$.
\endproclaim

\demo{Proof}
Let $A$ be a very ample line bundle on $T$.
If $\dim T > \dim f$ and $T_1$ is a general
member of $|A|$, then
$\rest{\phi}{T_1 \times Y} : T_1 \times Y \ratmap X$ still
dominates $X$. Thus, considering induction on $\dim T$,
we may assume that $\dim T = \dim f$.

Let $\mu : Z \to T \times Y$ be a birational morphism of
smooth projective varieties such that
$\psi = \phi \cdot \mu : Z \to X$ is a morphism.
Then, $\psi$ is generically finite. Thus, there is
a natural injection
$\psi^*(\omega_{X/Y}) \hookrightarrow \omega_{Z/Y}$.
Hence, $\psi^*(\omega_{X/Y}^n) \hookrightarrow \omega_{Z/Y}^n$
for all $n > 0$. Therefore,
$$
\omega_{X/Y}^n \hookrightarrow \psi_*(\psi^*(\omega_{X/Y}^n))
\hookrightarrow \psi_*(\omega_{Z/Y}^n).
$$
Applying $f_*$ to the above injection, we have
$$
f_*(\omega_{X/Y}^n) \hookrightarrow f_*(\psi_*(\omega_{Z/Y}^n)).
$$
Further,
$$
f_*(\psi_*(\omega_{Z/Y}^n)) = p_*(\mu_*(\omega_{Z/Y}^n))
= p_*(\omega_{T \times Y/Y}^n) = H^0(T, \omega_T^n) \otimes_k \OO_Y.
$$
Thus, $f_*(\omega_{X/Y}^n)^{\vee\vee}$ is a subsheaf of the free sheaf
$H^0(T, \omega_T^n) \otimes_k \OO_Y$.

Here we claim
$$
\left(c_1\left(f_*(\omega_{X/Y}^n)^{\vee\vee}\right)
\cdot H^{d-1}\right) \geq 0, \tag{1.1.1}
$$
where $H$ is an ample line bundle on $Y$ and $d = \dim Y$.
This is an immediate consequence of weak positivity of
$f_*(\omega_{X/Y}^n)^{\vee\vee}$ due to Viehweg \cite{Vi}.
We can however conclude our claim by a weaker result of Kawamata \cite{Ka1},
namely $\deg(f_*(\omega_{X/Y}^n)) \geq 0$ if $\dim Y = 1$.
For, considering complete intersections by general members of $|H^m|$ ($m \gg
0$),
we may assume $\dim Y = 1$.

We can find a projection $\alpha : H^0(T, \omega_T^n) \otimes_k \OO_Y
\to \OO_Y^{\oplus r_n}$ such that
$r_n = \rank f_*(\omega_{X/Y}^n)^{\vee\vee}$
and the composition
$$
f_*(\omega_{X/Y}^n)^{\vee\vee} \hookrightarrow
H^0(T, \omega_T^n) \otimes_k \OO_Y @>{\alpha}>> \OO_Y^{\oplus r_n}
$$
is injective. Therefore, since $f_*(\omega_{X/Y}^n)^{\vee\vee}$
is reflexive, the above homomorphism is an isomorphism by (1.1.1).
\QED
\enddemo

As corollary, we have a criterion of birational splitting.

\proclaim{Corollary 1.2}
Let $f : X \ratmap Y$ be a dominant rational map of algebraic varieties
over an algebraically closed field $k$ of characteristic zero,
whose generic fiber of $f$ is geometrically irreducible and of general type.
If there are an algebraic variety $T$ over $k$ and a dominant rational
map $\phi : T \times_k Y \ratmap X$ over $Y$, then there
is an algebraic variety $W$ over $k$ such that $X$ is birationally equivalent
to
a product $W \times_k Y$ over $Y$.
\endproclaim

\demo{Proof}
Clearly we may assume that $X$, $Y$ and $T$ are smooth and projective, and that
$f$ is a morphism. Moreover, in the same way as in the proof of Lemma~1.1,
we may assume that $\dim T = \dim f$.

Let $\mu : Z \to T \times Y$ be a birational morphism of
smooth projective varieties such that
$\psi = \phi \cdot \mu : Z \to X$ is a morphism.
We take a sufficiently large integer $m$ and an open set $U$ of $Y$ such that
\roster
\item "(a)" $|\omega_{f^{-1}(t)}^m|$ gives a birational map for all $t \in U$,
and

\item "(b)" the natural homomorphism
$f_*(\omega_{X/Y}^m) \otimes k(t) \to H^0(f^{-1}(t), \omega_{f^{-1}(t)}^m)$
is bijective for all $t \in U$, where $k(t)$ is the residue field at $t$.
\endroster
Let $U_0$ be an open set of $Y$ such that
$\codim_Y(Y \setminus U_0) \geq 2$ and
$$
\rest{f_*(\omega_{X/Y}^m)^{\vee\vee}}{U_0} =
\rest{f_*(\omega_{X/Y}^m)}{U_0}.
$$
We set $X_0 = f^{-1}(U_0)$.
Then we get
$$
H^0(X_0, \omega_{X/Y}^m) =
H^0\left(Y, f_*(\omega_{X/Y}^m)^{\vee\vee}\right).
\tag{1.2.1}
$$
Let $h : X_0 \ratmap \PP(H^0(X_0, \omega_{X/Y}^m))$ be a rational map
induced by $H^0(X_0, \omega_{X/Y}^m)$ and
$W$ the image of $h$. In order to see that
$h \times f : X_0 \ratmap W \times U_0$ is a birational map over $U_0$,
it is sufficient to show that $\rest{h}{f^{-1}(t)} : f^{-1}(t) \ratmap W$ is
a birational map for all $t \in U_0 \cap U$.

Here we claim $\dim W \leq \dim T$. First of all,
$$
\align
H^0(X_0, \omega_{X/Y}^m) & \hookrightarrow
H^0(\psi^{-1}(X_0), \omega_{Z/Y}^m) \\
& = H^0(T \times_k U_0, \omega_{T \times U_0 / U_0}^m) \\
& = H^0(T, \omega_T^m).
\endalign
$$
Thus if $h' : T \ratmap W$ is a rational map induced by
the image of
$H^0(X_0, \omega_{X/Y}^m) \hookrightarrow H^0(T, \omega_T^m)$,
then we have the following commutative diagram of rational maps:
$$
\CD
T \times_k U_0 @>{\phi}>> X_0 \\
@V{q}VV                 @VV{h}V \\
T            @>>{h'}>  W
\endCD
$$
where $q$ is the natural projection.
Since $\phi$, $h$ and $q$ are dominant rational maps,
so is $h'$. Therefore $\dim W \leq \dim T$.

By Lemma~1.1, $f_*(\omega_{X/Y}^m)^{\vee\vee}$
is a free sheaf on $Y$. Hence,
by virtue of (b) and (1.2.1),
$$
H^0(X_0, \omega_{X/Y}^m) \to H^0(f^{-1}(t), \omega_{f^{-1}(t)}^m)
$$
is bijective. This means that $\rest{h}{f^{-1}(t)} : f^{-1}(t) \ratmap W$ is
given by the complete linear system $|\omega_{f^{-1}(t)}^m|$.
Thus, by (a) and $\dim W \leq \dim f^{-1}(t)$,
$\rest{h}{f^{-1}(t)}$ is a birational map.
\QED
\enddemo

\head \S2. Proof of Proposition~B
\endhead

Clearly we may assume that $X$ and $Y$ are projective.
Let $\Rat_k(Y, X)$ be a scheme consisting of rational maps from $Y$ to $X$.
We have the natural morphism of schemes
$\alpha : \Rat_k(Y, X) \to \Rat_k(Y, Y)$
defined by $\alpha(s) = f \cdot s$
for $s \in \Rat_k(Y, X)$. We set
$\QSec(f) = \alpha^{-1}(\id_Y)$.
Then, it is easy to check the following properties of $\QSec(f)$.
\roster
\item "(a)" $\QSec(f)(k)$ is the set of all quasi-sections of $f : X \to Y$.

\item "(b)" There is the evaluation map
$\phi : \QSec(f) \times Y \ratmap X$,
which satisfies $f \cdot \phi = p$,
where $p : \QSec(f) \times Y \to Y$ is the natural projection.

\item "(c)" $\QSec(f)$ has only countably many connected components.
\endroster
(a) and (b) are trivial by our construction.
Concerning (c), $\Rat_k(Y, X)$ can be realized as an a open subscheme of
$\Hilb_{Y \times X}$. If we fix a polarization of $Y \times X$ and
a Hilbert polynomial $P$, then $\Hilb_{Y \times X}^P$ has only finite
connected components. On the other hand,
the Hilbert polynomials $P$ is a polynomials with coefficients in $\QQ$.
Thus we have only countably many possibilities of $P$.
Hence we can conclude (c).

Let $\QSec(f) = \bigcup_i \QSec_i(f)$ be the decomposition
into connected components.
Then, by the criterion of birational splitting, the evaluation map
$\QSec_i(f) \times Y \ratmap X$ is not a dominant rational map for every $i$.
Thus, we get our proposition.
\QED

\head
\S3. Northcott's theorem over function fields
\endhead

By the same idea as the proof of Proposition~B,
we have the following Northcott's theorem over
function fields, which is a generalization of Theorem~1.3 of \cite{Mo}.

\proclaim{Proposition 3.1}
Let $X$ be a smooth projective variety over an algebraically closed field $k$
of characteristic zero,
$C$ a smooth projective curve over $k$, and $f : X \to C$
a surjective morphism whose generic fiber is geometrically irreducible.
Let $L$ be a line bundle on $X$ such that
$L$ is ample on the generic fiber of $f$.
Assume that
$\deg(f_*(\omega_{X/C}^{n})) >0$ for some $n > 0$.
Then, for any number $A$, the set
$$
\left\{ \Delta \mid \hbox{$\Delta$ is a section of $f : X \to C$ with
$(L \cdot \Delta) \leq A$}
\right\}
$$
is not dense in $X$.
\endproclaim

\demo{Proof}
In the same way as in the proof of Theorem~1.3 of \cite{Mo},
we may assume that $L$ is ample on $X$.
Thus, the set
$$
S = \left\{ \Delta \mid \hbox{$\Delta$ is a section of $f : X \to C$ with
$(L \cdot \Delta) \leq A$}
\right\}
$$
is a bounded family. Therefore, there is a finite union $Y$ of
irreducible components of $\QSec(f)$ such that
$Y(k)$ coincides with $S$.
By Lemma~1.1, the evaluation map
$\rest{\phi}{Y \times C} : Y \times C \to X$ doesn't dominate $X$.
Thus, we get our proposition.
\QED
\enddemo

\head \S4 Remarks
\endhead
\subhead Remark 4.1
\endsubhead
Let $k$ be an algebraically closed field of characteristic zero.
Let $f : X \to Y$ be a projective and surjective morphism
of algebraic varieties over $k$.
Let $K$ be the algebraic closure of the function field $k(Y)$ of $Y$.
The variation of $f$, denoted by $\Var(f)$, is defined by
the minimal transcendental degree of a field $L$ such that
$k \subset L \subset K$ and $X_{K}$
is birationally equivalent to $W_{K}$ for some
projective variety $W$ over $L$.
The fundamental conjecture of classification theory of algebraic varieties
is the following:
\block
Let $f : X \to Y$ be a morphism of smooth projective varieties
over $k$ with connected fibers. If the Kodaira dimension of
the generic fiber is non-negative, then, for sufficiently large $n$,
$\kappa(Y, \det(f_*(\omega_{X/Y}^n))^{\vee\vee}) \geq \Var(f)$.
\endblock
For example, the above is known if the geometric generic fiber has a good
minimal
model (cf. \cite{Ka2}).
You can find more details in \cite{M}.
Using Lemma~2.1, we can see that the above conjecture implies
the following:
\block
Let $f : X \to Y$ be a projective and surjective morphism of algebraic
varieties
over $k$, whose generic fiber is geometrically irreducible.
If the generic fiber of $f$ has non-negative Kodaira dimension and
the variation of $f$ is positive,
then there are countably many proper closed varieties
$\{ Z_i \}$ of $X$ such that
every quasi-section of $f$ is contained in $\bigcup_i Z_i$.
\endblock
Therefore, we can make a generalization of Conjecture~A, namely,
\block
Let $f : X \to Y$ be a projective and surjective morphism of algebraic
varieties
over $k$, whose generic fiber is geometrically irreducible.
If the generic fiber of $f$ has non-negative Kodaira dimension and
the variation of $f$ is positive,
then there are a proper subscheme $Z$ of $X$
such that every quasi-section of $f$ is contained in $Z$.
\endblock

\subhead Remark 4.2
\endsubhead
Let $k$ be an algebraically closed field of characteristic $p > 0$.
In \cite{Sh}, T. Shiota constructed a family of unirational hypersurfaces in
$\PP_k^3$.
More precisely, letting $q = p^{\nu}$ ($\nu \geq 1$),
he considered surfaces of degree $q+1$ in $\PP_k^3$ defined by
$$
x_1^qx_3 + x_2^qx_4 + x_1f(x_3, x_4) + x_2g(x_3, x_4) = 0,
$$
where $f$ and $g$ are binary forms of degree $q$ in $x_3$ and $x_4$
without common factors.
He checked they are unirational by direct calculations.
He also checked the number of essential parameters of the above type surfaces
is $2q-2$.
Since they have ample canonical line bundles,
a birational map between them is an isomorphism.
Therefore this family is not birationally trivial.
This observation shows us that Conjecture~A does not hold in this example.


\enddocument